\documentclass[12pt,leqno]{article}
\title{On  $\pi$-adic expansion  of singular integers  of the  $p$-cyclotomic field}
\author{Roland Qu\^eme}
\usepackage{amsmath,amsthm,amstext,amsbsy}
\newtheorem{thm}{Theorem}[section]

\newtheorem{lem}[thm]{Lemma}
\font\mathbb=msbm10
\newcommand{\N}{\mbox{\mathbb N}}

\newcommand{\Q}{\mbox{\mathbb Q}}
\newcommand{\Z}{\mbox{\mathbb Z}}

\newcommand{\modu}{\ \mbox{mod}\ }
\newcommand{\be}{\begin{equation}}
\newcommand{\ee}{\end{equation}}
\date{2006 nov 12}
\begin{document}
\maketitle
%\tableofcontents \clearpage
\abstract
Let $p$ be an odd prime. Let $K=\Q(\zeta)$ be the $p$-cyclotomic field and let $O_K$ be the ring of integers of $K$.
Let $\pi$ be the prime ideal of $K$ lying over $p$.
An integer  $B\in O_K$ is said singular if $B^{1/p}\not\in K$ and if $B O_K =\mathbf b^p$ where $\mathbf b$ is an  ideal of $O_K$.
An integer  $B\in O_K$ is said semi-primary if $B\equiv b\modu\pi^2$ where $b\in \Z,\ b\not\equiv 0\modu p$.
Let $\sigma$ be a $\Q$-isomorphism of the field $K$ generating the Galois group $Gal(K/\Q)$.
When  $p$ is irregular, there exists at least one subgroup  $\Gamma$  of order $p$ of the class group of $K$ annihilated by  a polynomial
$\sigma-\mu$ with  $\mu\in {\bf F}_p^*$.
We prove   the existence, for each $\Gamma$,  of singular semi-primary integers $B$ where $B O_K=\mathbf b^p$ with class $Cl(\mathbf b)\in \Gamma$
and  $B^{\sigma-\mu}\in K^p$ and we describe their $\pi$-adic expansion.
This paper is at a strictly elementary level.
\endabstract
%%% ====================================================================
%
%RRRRRRRR 20
%%% ====================================================================
%%% ====================================================================
%
%RRRRRRRR 20
%%% ====================================================================
\section{Some definitions on cyclotomic fields}
In this section, we fix some definitions and notations and remind reader  of some classical properties of cyclotomic fields used in  the article.
\begin{enumerate}
\item
Let $p$ be an odd prime.
Let  $K=\Q(\zeta)$ be the $p$-cyclotomic number field.
Let   $O_K=\Z[\zeta]$ be the ring of integers of $K$.
Let $K^+=\Q(\zeta+\zeta^{-1})$ be the maximal totally real subfield of $K$.
The ring of integers of $K^+$ is $O_{K^+}=\Z[\zeta+\zeta^{-1}]$. Let us denote $O_{K^+}^*$ the group of units of
$O_{K^+}$.
Let ${\bf F}_p$ be the finite field with $p$ elements and ${\bf F}_p^*$ its multiplicative group.
\item
Let us denote $\mathbf a$ the  ideals of $O_K$ and  $Cl(\mathbf a)$ their  classes in the class group of $K$.
Let us denote $<Cl(\mathbf a)>$  the finite group generated by the class $Cl(\mathbf a)$.
If $a\in O_K$ then  $a O_K$ is the  principal  ideal of $O_K$ generated by $a$.
The ideal $p O_K=\pi^{p-1}$ where $\pi$ is the principal prime ideal $(1-\zeta)O_K$.
Let us denote $\lambda= \zeta-1$, so $\pi=\lambda O_K$.
\item
Let $G=Gal(K/\Q)$ be the Galois group of the field $K$.
Let $\sigma$ be a $K$-isomorphism generating the cyclic group $G$.
$\sigma$ is  defined by $\sigma(\zeta)=\zeta^u$ where $u$ is a primitive root $\modu p$.
\item
For this primitive root $u\modu p$ and $i\in\N$,
let us denote $u_i\equiv u^i \modu p,\quad 1\leq u_i\leq p-1$.
For $i\in\Z,\quad i<0$, this is to be understood as $u_i u^{-i}\equiv 1 \modu p$.
This notation follows the convention adopted in Ribenboim \cite{rib}, last paragraph of page 118.
This notation is largely used in the sequel of this article.
\item
Let $C_p, C_p^+$ be  the subgroups of exponent $p$ of the class groups of the field $K$ and $K^+$.
Then $C_p=C_p^+\oplus C_p^-$ where  $C_p^-$ is the relative $p$-class group $C_p/C_p^+$.
Let $r$ be   the rank of the groups  $C_p$.
The abelian group $C_p$ is a group of order $p^{r}$ with
\begin{equation}\label{e611084}
C_p=\oplus_{i=1}^{r} \Gamma_i,
\end{equation}
where each $\Gamma_i$ is a cyclic groups of order $p$ annihilated by $\sigma-\mu_i$ with $\mu_i\in {\bf F}_p^*$ .
\item
Similarly   the $p$-unit group $U=O_{K^+}^*/O_{K^+}^{*p}$ is a direct sum
\begin{equation}\label{e611085}
U=\oplus _{i=1}^{(p-3)/2} U_i,
\end{equation}
where $U_i=<\eta_i>$ is a cyclic group of order $p$ with
$\eta_i^{\sigma-\mu_i}\in O_{K^+}^{*p}$ and  where $\mu_i= u_{2m_i}\modu p$ with $1 \leq m_i\leq\frac{p-3}{4}$.
\item
We say that an  algebraic number  $C\in K$ is  singular if  $C^{1/p}\not\in K$ and $C O_K =\mathbf c^p$ for
some ideal $\mathbf c$ of $K$. If $C$ is integer then $C$ is called a singular integer.
Observe that with this definition a unit is a singular integer.
We say that $C$ is singular primary
if $C$ is singular and $C\equiv c^p\modu \pi^p,\quad c\in \Z,\quad c\not\equiv 0\modu p$.
We say that  the singular number $C$  is   semi-primary if $C\equiv c\modu\pi^2,\quad c\in\Z,\quad c\not\equiv 0\modu p$.
Observe that if $C$ is primary then $C$ is semi-primary.
\item
When  $p$ is irregular, there exists at least one subgroup  $\Gamma$  of order $p$ of the class group of $K$ annihilated by  a polynomial
$\sigma-\mu$ with  $\mu\in {\bf F}_p^*$.
In this article we prove, for each $\Gamma$, the existence of singular semi-primary integers $B$ where $B O_K=\mathbf b^p$ with class
$Cl(\mathbf b)\in \Gamma$
and  $B^{\sigma-\mu}\in K^p$ and we describe their $\pi$-adic expansion.
\end{enumerate}
%%% ====================================================================
%
%RRRRRRRR 260
%%% ====================================================================
%%% ====================================================================
%
%RRRRRRRR 260
%%% ====================================================================
%%% ====================================================================
%
%RRRRRRRR 260
%%% ====================================================================
\section{On $\pi$-adic expansion of singular integers}\label{s108311}
In this section we consider  the singular integers  $B$ with   $B O_K=\mathbf b^p$ where  $Cl(\mathbf b)$ is annihilated by $\sigma-\mu$ for one
$\mu\in{\bf F}_p^*$. The singular integer $B$ is said respectively {\it negative}  when  $Cl(\mathbf b)\in C_p^-$ and {\it positive}
when $Cl(\mathbf b)\in C_p^+$. Observe that $\mu^{(p-1)/2}\equiv -1\modu p$ when $Cl(\mathbf b)\in C_p^-$ and that
 $\mu^{(p-1)/2}\equiv 1\modu p$ when $Cl(\mathbf b)\in C_p^+$.
%%% ====================================================================
%
%RRRRRRRR 260
%%% ====================================================================
\subsection{$\pi$-adic expansion of singular negative integers}
%%% ====================================================================
%
%RRRRRRRR 260
%%% ====================================================================
At first,  we give a general lemma dealing with congruences on $p$-powers of algebraic numbers of $K$.
\begin{lem}$ $\label{l502031}
Let $\alpha,\beta\in O_K$
with $\alpha\not\equiv 0 \mod \pi$
and $\alpha\equiv\beta\modu \pi$.
Then $\alpha^p\equiv \beta^p  \modu \pi^{p+1}$.
\begin{proof}$ $
Let $\lambda=(\zeta-1)$.
Then
$\alpha-\beta\equiv 0 \modu \pi$  implies that
$\alpha-\zeta^k\beta\equiv 0 \modu \pi$ for $k=0,1,\dots,p-1$.
Therefore, for all
$k, \quad 0\leq k \leq p-1$, there exists $a_k\in \N,\quad 0\leq a_k \leq p-1$, such that
$(\alpha-\zeta^k\beta)\equiv \lambda a_k \modu \pi^2$.
For another value $l,\quad 0\leq l\leq p-1$, we have, in the same way,
$(\alpha-\zeta^l\beta)\equiv \lambda a_l \modu \pi^2$,
hence $(\zeta^k-\zeta^l)\beta\equiv \lambda(a_k-a_l) \modu \pi^2$.
For $k\not=l$ we get $a_k\not= a_l$, because $\pi \| (\zeta^k-\zeta^l)$ and
because hypothesis $\alpha\not\equiv 0 \modu \pi$ implies that $\beta\not\equiv 0 \modu \pi$.
Therefore, there exists one and only one $k$ such that
$(\alpha-\zeta^k\beta)\equiv 0 \modu \pi^2$. Then, we have
$ \prod_{j=0}^{p-1}(\alpha-\zeta^j\beta)
= (\alpha^p-\beta^p)\equiv 0 \modu \pi^{p+1}$.
\end{proof}
\end{lem}
%%% ====================================================================
%
%RRRRRRRR 30
%%% ====================================================================
\begin{lem}\label{l108161}
Let $\mathbf b$ be an ideal of $O_K$ such that $Cl(\mathbf b)\in C_p^-$ is annihilated by $\sigma-\mu$.
There exist    singular semi-primary negative integers $B$  with $B O_K=\mathbf b^p$. They  verify the relation
\begin{equation}\label{e610125}
\begin{split}
& (\frac{B}{\overline{B}})^{\sigma-\mu}
=(\frac{\alpha}{\overline{\alpha}})^p,\quad\alpha\in K.\\
\end{split}
\end{equation}
\begin{proof}$ $
The ideal $\mathbf b^p$ is principal. So let one $\beta\in\Z[\zeta]$ with
$\beta O_K=\mathbf b^p$. There exists a  natural number $w$ such that the integer $B=\beta\zeta^w$ is semi-primary.
$\mathbf b^{\sigma-\mu}$ is principal, therefore there exists
$\alpha\in K$ such that
$\sigma(\mathbf b) =\mathbf b^{\mu}\alpha O_K$,
thus there exists a natural number $w^\prime$ with
\begin{displaymath}
\sigma( B) = B^{\mu} \eta\zeta^{w^\prime} \alpha^p,\quad \eta\in O_{K^+}^*, \quad w^\prime\in \Z.
\end{displaymath}
$B$ is semi-primary, hence $\sigma(B)$ is semi-primary. $\eta$ and $\alpha^p$ are semi-primary, hence $w^\prime=0$.
Then $\sigma( B) = B^{\mu} \eta \alpha^p$ and by conjugation
$\sigma(\overline{ B}) =\overline{ B}^{\mu} \eta \overline{\alpha}^p$ and the result follows.
\end{proof}
\end{lem}
%%% ====================================================================
%
%RRRRRRRR 260

%%% ====================================================================
\begin{lem}\label{l108171}
The singular semi-primary negative integer $B$ defined in previous lemma \ref{l108161}  verifies $\mu=u_{2m+1}$
for a natural integer $\ m,\quad  1\leq m\leq\frac{p-3}{2}$.
If  the singular number $C=\frac{B}{\overline{B}}$ is non-primary then   $\pi^{2m+1}\ \|\ C-1$.
\begin{proof} $ $
\begin{enumerate}
\item
$Cl(\mathbf b)\in C_p^-$, hence $\mu^{(p-1)/2}\equiv -1\modu p$, thus $\mu=u_{2m+1}$ for a  natural number
$m,\ 0\leq m\leq \frac{p-3}{2}$. From Stickelberger theorem $\sigma-u$ does not annihilate $Cl(\mathbf b)$, hence
$m\not=0$.
\item
The definition of $C$ implies that $C\equiv 1\modu \pi$.
Observe that $v_\pi(C-1)$ is odd because $C\overline{C}=1$. Therefore the hypothesis $C$ non-primary implies that $v_\pi(C-1)\leq p-2$.
There exists a natural integer $\nu$ such that $\pi^\nu\ \|\ C-1$, hence
\begin{equation}\label{e502031}
\begin{split}
& C\equiv 1+c_0\lambda^\nu\modu \lambda^{\nu+1},\ \nu\leq p-2,\\
& c_0\in\Z,\quad c_0\not\equiv 0\modu p.\\
\end{split}
\end{equation}
We have to prove that $\nu\leq p-2$ implies that
$\nu=2m+1$:
\item
From lemma \ref{l108161} p.\pageref{l108161}, it follows that $\sigma(C)=C^\mu\times\alpha_1^p$, with
$\alpha_1=\frac{\alpha}{\overline{\alpha}}$, and so that
$1+c_0\sigma(\lambda)^\nu\equiv  (1+\mu c_0\lambda^\nu)\times \alpha_1^p\modu\pi^{\nu+1}$.
In the other hand $\alpha_1\equiv 1\modu\pi$ and then, from lemma \ref{l502031},  $\alpha_1^p\equiv 1\modu\pi^{p+1}$.
Then $1+c_0\sigma(\lambda)^\nu\equiv 1+\mu c_0\lambda^\nu\modu\lambda^{\nu+1}$, and so
$\sigma(\lambda^{\nu})\equiv \mu\lambda^{\nu}\modu\pi^{\nu+1}$. This implies that $\sigma(\zeta-1)^\nu\equiv \mu\lambda^{\nu}\modu\pi^{\nu+1}$,
so that $(\zeta^u-1)^\nu\equiv \mu\lambda^{\nu}\modu\pi^{\nu+1}$, so that
$((\lambda+1)^u-1)^\nu\equiv \mu\lambda^{\nu}\modu\pi^{\nu+1}$ and finally
$u^\nu\lambda^\nu\equiv \mu\lambda^{\nu}\modu\pi^{\nu+1}$, hence  $u^\nu-\mu\equiv 0\modu \pi$.
Therefore, we have proved that $\nu=2m+1$.
\end{enumerate}
\end{proof}
\end{lem}
%%% ====================================================================
%
%RRRRRRRR 260
%%% ====================================================================

%%% ====================================================================
%
%RRRRRRRR 260
%%% ====================================================================
In this theorem we generalize to singular integers $B$ the results obtained for singular numbers $C$ in lemma \ref{l108161} and \ref{l108171}.
\begin{thm} \label{l203171}$ $
Let $B$ be a singular semi-primary  negative integer.
There exist    a unit $\eta\in O_{K^+}^*$ given by $B\overline{B}=\eta\times\beta^p,\ \beta\in O_K$,  and a singular  negative integer
$B^{\prime}= \frac{B^2}{\eta}$     such that
\begin{equation}\label{e203271}
\begin{split}
& \sigma(B^{\prime})=B^{\prime\mu}\times \alpha^{\prime p},\quad \alpha^\prime\in K,\\
\end{split}
\end{equation}
If $B^\prime$ is non-primary then $\pi^{2m+1}\ \|\ (B^{\prime})^{p-1}-1$.
\begin{proof}$ $
From lemma \ref{l108161} we have $(\frac{B}{\overline{B}})^{\sigma-\mu}=(\frac{\alpha}{\overline{\alpha}})^p$.
$Cl(\mathbf b)\in C_p^-$ implies that $\mathbf b\overline{\mathbf b}$ is principal
and so that $B\overline{B}=\eta\beta^p$ with $\beta\in O_K$ and $\eta\in O_{K^+}^*$.
Therefore
\begin{displaymath}
(\frac{B^2}{\eta})^{\sigma-\mu}=(\frac{B^2}{\frac{B\overline{B}}{\beta^p}})^{\sigma-\mu}=(\frac{B}{\overline{B}}\times \beta^p)^{\sigma-\mu}
=(\frac{\alpha}{\overline{\alpha}}\times \beta^{\sigma-\mu})^p.
\end{displaymath}
Let us denote $B^\prime=\frac{B^2}{\eta},\quad B^\prime\in O_K,\quad v_\pi(B^\prime)=0$.
We get
\begin{equation}\label{e203092}
\sigma(B^\prime)=(B^\prime)^\mu\times \alpha^{\prime p},\  \alpha^\prime\in K.
\end{equation}
This relation leads to
\begin{equation}\label{e610251}
\sigma(B^\prime )^{p-1}\equiv (B^\prime)^{(p-1)\mu}\modu\pi^{p+1}.
\end{equation}
If $B^\prime$ is non-primary then it leads in the same way  than in
lemma  \ref{l108171} p. \pageref{l108171} to  the congruence
\begin{displaymath}
\pi^{2m+1}\ \|\ (B^{\prime})^{p-1}-1,
\end{displaymath}
 which achieves  the proof.
\end{proof}
\end{thm}
\paragraph{Remark:} If  $B$ is singular primary then $\eta=1$ from a  theorem of Furtwangler, see Ribenboim \cite{rib} (6C) p. 182.
%%% ====================================================================
%
%RRRRRRRR 260
%%% ====================================================================
%%% ====================================================================
%
%RRRRRRRR 260
%%% ====================================================================
\paragraph {The case $\mu=u_{2m+1}$ with $2m+1>\frac{p-1}{2}$}\label{s201111}
\paragraph{ }
Let us consider the singular number $C$ defined in lemma \ref{l108171} p. \pageref{l108171}. The number $C$ can be written in the form
\begin{displaymath}
\begin{split}
& C=1+\gamma+\gamma_0\zeta+\gamma_1 \zeta^u
+\dots + \gamma_{p-3}\zeta^{u_{p-3}},\\
& \gamma\in\Q,\quad v_p(\gamma)\geq 0,\quad  \gamma_i\in\Q,\quad
v_p(\gamma_i)\geq 0,\quad i=0, \dots, p-3,\\
& \gamma+\gamma_0\zeta+\gamma_1\zeta^u+\dots + \gamma_{p-3}\zeta^{u_{p-3}}\equiv 0\modu \pi^{2m+1}.
\end{split}
\end{displaymath}
In the   theorem \ref{l202211} p. \pageref{l202211},
we shall compute    the coefficients  $\gamma$ and $\gamma_i\modu p$.
%%% ====================================================================
%
%RRRRRRRR 260
%%% ====================================================================
\begin{lem}\label{l201091}
$C$ verifies  the congruences
\begin{displaymath}
\begin{split}
& \gamma\equiv -\frac{\gamma_{p-3}}{\mu-1}\modu p,\\
& \gamma_0\equiv -\mu^{-1} \times \gamma_{p-3}\modu p,\\
& \gamma_1\equiv -(\mu^{-2}+\mu^{-1})\times \gamma_{p-3}\modu p,\\
&\vdots\\
& \gamma_{p-4}\equiv -(\mu^{-(p-3)}+\dots+\mu^{-1})\times \gamma_{p-3}\modu p.
\end{split}
\end{displaymath}
\begin{proof}
We have seen in lemma \ref{l108161} p. \pageref{l108161} that $\sigma(C)\equiv C^{\mu}\modu\pi^{p+1}$.
From $2m+1 > \frac{p-1}{2}$ we derive that
\begin{displaymath}
C^\mu\equiv 1+\mu\times(\gamma+\gamma_0\zeta+\gamma_1\zeta^u
+\dots +\gamma_{p-3} \zeta^{u_{p-3}})\modu\pi^{p-1}.
\end{displaymath}
In the other hand, we get by conjugation
\begin{equation}\label{e201111}
\sigma(C)=1+\gamma+\gamma_0\zeta^u+\gamma_1\zeta^{u_2}
+\dots +\gamma_{p-3} \zeta^{u_{p-2}}.
\end{equation}
We have the identity
\begin{displaymath}
\gamma_{p-3} \zeta^{u_{p-2}}=-\gamma_{p-3}-\gamma_{p-3}\zeta
-\dots -\gamma_{p-3}\zeta^{u_{p-3}}.
\end{displaymath}
This leads to
\begin{displaymath}
\sigma(C)=1+\gamma-\gamma_{p-3}-\gamma_{p-3}
\zeta+(\gamma_0-\gamma_{p-3})\zeta^u
+\dots+(\gamma_{p-4}-\gamma_{p-3})\zeta^{u_{p-3}}.
\end{displaymath}
Therefore, from the congruence $\sigma(C)\equiv C^{\mu}\modu \pi^{p+1}$ we get
the congruences in the basis $1,\zeta,\zeta^{u},\dots,\zeta^{u_{p-3}}$,
\begin{displaymath}
\begin{split}
& 1+\mu \gamma\equiv 1+\gamma-\gamma_{p-3}\modu p,\\
& \mu \gamma_0\equiv -\gamma_{p-3}\modu p,\\
& \mu \gamma_1\equiv \gamma_0-\gamma_{p-3}\modu p,\\
& \mu \gamma_2\equiv \gamma_1-\gamma_{p-3}\modu p,\\
& \vdots\\
& \mu \gamma_{p-4}\equiv \gamma_{p-5}-\gamma_{p-3}\modu p,\\
& \mu \gamma_{p-3}\equiv \gamma_{p-4}-\gamma_{p-3}\modu p.
\end{split}
\end{displaymath}
From these congruences, we get $\gamma\equiv -\frac{\gamma_{p-3}}{\mu-1}\modu p$ and
$\gamma_0\equiv -\mu^{-1} \gamma_{p-3}\modu p$ and then
$\gamma_1\equiv \mu^{-1}(\gamma_0-\gamma_{p-3})\equiv \mu^{-1}(-\mu^{-1}
\gamma_{p-3}-\gamma_{p-3})
\equiv -(\mu^{-2}+\mu^{-1})\gamma_{p-3}\modu p$
and $\gamma_2
\equiv\mu^{-1}(\gamma_1-\gamma_{p-3})
\equiv \mu^{-1}(-(\mu^{-2}+\mu^{-1})\gamma_{p-3}-\gamma_{p-3})
\equiv -(\mu^{-3}+\mu^{-2}+\mu^{-1})\gamma_{p-3}\modu p$ and so on.
\end{proof}
\end{lem}
%%% ====================================================================
%
%RRRRRRRR 260
%%% ====================================================================
\begin{thm}{}\label{l202211}
If $2m+1> \frac{p-1}{2}$ then $C$ verifies the congruence
\begin{equation}\label{e202211}
C\equiv 1-\delta \times
(\zeta+\mu^{-1}\zeta^u+\dots+\mu^{-(p-2)}\zeta^{u_{p-2}}) \modu \pi^{p-1},
\end{equation}
where $\delta\in\Z$ is coprime with $p$ when $C$ is non-primary.
\begin{proof}
The result is trivial if $C$ is primary. Suppose that $C$ is not primary.
From definition of $C$, setting $C=1+V$, we get :
\begin{displaymath}
\begin{split}
& C=1+V,\\
& V=\gamma+\gamma_0\zeta+\gamma_1\zeta^u+\dots+\gamma_{p-3}\zeta^{u_{p-3}},\\
& \sigma(V)\equiv \mu\times V\modu \pi^{p+1}.
\end{split}
\end{displaymath}
Then, from lemma \ref{l201091} p. \pageref{l201091},  observing that $\mu^{-(p-2)}+\dots+\mu^{-1}\equiv -1\modu p$, we obtain the relations
\begin{displaymath}
\begin{split}
& \mu=u_{2m+1},\\
& \gamma \equiv -\frac{\gamma_{p-3}}{\mu-1}\modu p,\\
& \gamma_0\equiv -\mu^{-1} \times \gamma_{p-3}\modu  p,\\
& \gamma_1\equiv -(\mu^{-2}+\mu^{-1})\times \gamma_{p-3}\modu p,\\
&\vdots\\
& \gamma_{p-4}\equiv -(\mu^{-(p-3)}+\dots+\mu^{-1})\times \gamma_{p-3}\modu  p,\\
& \gamma_{p-3}\equiv -(\mu^{-(p-2)}+\dots+\mu^{-1})\times \gamma_{p-3}\modu  p.
\end{split}
\end{displaymath}
From these relations we get
\begin{displaymath}
V\equiv  -\gamma_{p-3}\times
(\frac{1}{\mu-1}+\mu^{-1}\zeta+(\mu^{-2}+\mu^{-1})\zeta^{u}+\dots+
(\mu^{-(p-2)}+\dots+\mu^{-1})\zeta^{u_{p-3}})\modu p,
\end{displaymath}
hence
\begin{displaymath}
V\equiv-\gamma_{p-3}\times (\frac{1}{\mu-1}+
\mu^{-1}(\zeta+(\mu^{-1}+1)\zeta^u+\dots+(\mu^{-(p-3)}+\dots+1)
\zeta^{u_{p-3}}))\modu p,
\end{displaymath}
hence
\begin{displaymath}
V\equiv -\gamma_{p-3}\times (\frac{1}{\mu-1}+
\mu^{-1}(
\frac{(\mu^{-1}-1)\zeta+(\mu^{-2}-1)\zeta^u+\dots
+(\mu^{-(p-2)}-1)\zeta^{u_{p-3}}}{\mu^{-1}-1}))
\modu p,
\end{displaymath}
hence
\begin{displaymath}
V\equiv -\gamma_{p-3}\times (\frac{1}{\mu-1}+
\mu^{-1}(
\frac{\mu^{-1}\zeta+\mu^{-2}\zeta^u+\dots
+\mu^{-(p-2)}\zeta^{u_{p-3}}
-\zeta-\zeta^u-\dots -\zeta^{u_{p-3}}}{\mu^{-1}-1}))
\modu p.
\end{displaymath}
In the other hand   $-\zeta-\zeta^u-\dots-\zeta^{u_{p-3}}=1+\zeta^{u_{p-2}}$ and $\mu^{-(p-1)}\equiv 1\modu p$
implies that
\begin{displaymath}
V\equiv -\gamma_{p-3}\times (\frac{1}{\mu-1}+
\mu^{-1}(
\frac{1+\mu^{-1}\zeta+\mu^{-2}\zeta^u+\dots
+\mu^{-(p-2)}\zeta^{u_{p-3}}+\mu^{-(p-1)}\zeta^{u_{p-2}}}{\mu^{-1}-1}))
\modu p,
\end{displaymath}
hence
\begin{displaymath}
V\equiv -\gamma_{p-3}\times (\frac{1}{\mu-1})\times
(1- (1+\mu^{-1}\zeta+\mu^{-2}\zeta^u+\dots
+\mu^{-(p-2)}\zeta^{u_{p-3}}+\mu^{-(p-1)}\zeta^{u_{p-2}}))
\modu p,
\end{displaymath}
hence
\begin{displaymath}
V\equiv -\gamma_{p-3}\times (\frac{\mu^{-1}}{\mu-1})\times
(\zeta+\mu^{-1}\zeta^u+\dots
+\mu^{-(p-3)}\zeta^{u_{p-3}}+\mu^{-(p-2)}\zeta^{u_{p-2}})
\modu p,
\end{displaymath}
which achieves the proof.
\end{proof}
\end{thm}
\paragraph{Remark:} we have a similar result with $(B^\prime)^{p-1}$ in place of $C$.
%%% ====================================================================
%
%RRRRRRRR 260
%%% ====================================================================
\subsection{$\pi$-adic expansion of singular positive integers}
Let $\mathbf b$ be an  ideal of $O_K$ whose class $Cl(\mathbf q)\in C_p^+$ is annihilated by
$\sigma-\mu$. In that case $\mu^{(p-1)/2}\equiv 1\modu p$ and $\mu = u_{2m}\modu p,\ 1\leq m\leq\frac{p-3}{2}$.
%%% ====================================================================
%
%RRRRRRRR 260
%%% ====================================================================
\begin{thm}  \label{l203272}$ $
\begin{enumerate}
\item
There exists singular semi-primary positive integers $B\in O_K$ such that:
\begin{equation}\label{e205041}
\begin{split}
& B O_K=\mathbf b^p,\\
& \sigma(B)=B^{\mu}\times\alpha^p,\quad\alpha\in K,\\
\end{split}
\end{equation}
\item
If $B$ is non-primary then $\pi^{2m}\ \|\ B^{p-1}-1$.
\end{enumerate}
\begin{proof}
There exists semi-primary integers $B^\prime$ with $B^\prime O_K=\mathbf b^p$ such that
\begin{displaymath}
\sigma(B^\prime)=B^{\prime\mu}\times\alpha^p\times \eta,
\quad \alpha\in K,\quad \eta\in O_{K^+}^*.
\end{displaymath}
From  independent   forward theorem \ref{t610122} p. \pageref{t610122}  dealing with
unit group $O_{K^+}^*$,  the unit $\eta$ verifies  the relation
\begin{equation}\label{e610121}
\begin{split}
& \eta= \eta_1^{l_1}\times (\prod_{j=2}^{N} \eta_{j}^{l^j}),
\quad l_j\in {\bf F}_p,\quad 1\leq  N<\frac{p-3}{2},\\
& \sigma(\eta_1)=\eta_1^{\mu}\times\beta_1^p,
\quad \eta_1,\beta_1\in O_{K^+}^*,\\
& \sigma(\eta_{j})=\eta_{j}^{\nu_{j}}\times\beta_j^p,
\quad \eta_j,\beta_j\in O_{K^+}^*, \quad j=2,\dots,N,\\
\end{split}
\end{equation}
where $\nu_j\not=\nu_{j^\prime}$ for  $2\leq j<j^\prime\leq N$ and $\nu_j\not=\mu$ for $j=2,\dots,N$.
Let us denote  $E=\eta_1^{l_1}$ and $U=\prod_{j=2}^{N} \eta_j^{l_j}$, hence $\eta=EU$.
Show that there exists $V\in O_{K^+}^*$ of form  $V=\prod_{j=2}^{N} \eta_j^{\rho_j}$ such that
\begin{equation}\label{e610122}
\sigma(V)\times V^{-\mu}=U^{-1}\times \varepsilon^p,\quad \varepsilon\in O_{K^+}^*
\end{equation}
It is sufficient   that
\begin{displaymath}
\eta_j^{\rho_j \nu_j}\times \eta_j^{-\rho_j\mu}
=\eta_j^{-l_j}\times \varepsilon_j^p,
\quad \varepsilon_j\in O_{K^+}^*,\quad j=2,\dots,N,
\end{displaymath}
hence   that
\begin{displaymath}
\rho_j\equiv \frac{-l_j}{\nu_j-\mu}\modu p,\quad j=2,\dots,N,
\end{displaymath}
which is possible, because $\nu_j\not\equiv \mu,\quad j=2,\dots, N$.
Therefore, for $B=B^\prime\times V$, we get  $B^\prime=B V^{-1}$ and so
\begin{displaymath}
\sigma(B^\prime)=\sigma(B V^{-1})=B^{\prime\mu}\alpha^p\eta=(B V^{-1})^\mu\alpha^p\eta=(B V^{-1})^{\mu}\times\alpha^p\times E\times U,
\end{displaymath}
hence
\begin{displaymath}
\sigma(B )=B^{\mu}\sigma(V) V^{-\mu}\times\alpha^p\times E\times U.
\end{displaymath}
From relation (\ref{e610122}) we get
\begin{displaymath}
\sigma(B )=B^{\mu}(U^{-1}\varepsilon^p\times\alpha^p\times E\times U),
\end{displaymath}
hence  we get the two simultaneous relations
\begin{equation}\label{e610253}
\begin{split}
& \sigma(B)=B^{\mu}\times \alpha^p\times\varepsilon^p\times E,\quad \alpha\in K,\\
& \sigma(E)=E^{\mu}\times \varepsilon_1^p,\quad \varepsilon_1\in O_{K^+}*.\\
\end{split}
\end{equation}
Show that
\begin{equation}\label{e407214}
B^{\sigma-\mu}\in K^{p}.
\end{equation}
\begin{enumerate}
\item
If $E\in  O_{K^+}^{*p}$, it is clear from  relation (\ref{e610253}).
\item
If $E\not\in  O_{K^+}^{*p}$
\begin{equation}\label{e611081}
\begin{split}
& \sigma(B)=B^{\mu}\times \alpha_1^{ p}\times E,\quad \alpha_1\in K,\\
\end{split}
\end{equation}
hence raising (\ref{e611081}) to $\mu$-power we get
\begin{equation}\label{e611082}
\begin{split}
& \sigma(B)^\mu=B^{\mu^2}\times \alpha_1^{ p\mu}\times E^\mu,\quad \alpha_1\in K,\\
\end{split}
\end{equation}
and also applying $\sigma$ to relation (\ref{e611081}) we get
\begin{equation}\label{e611083}
\begin{split}
& \sigma^2(B)=\sigma(B)^{\mu}\times  E^{\mu}\times b^p.
\quad b\in K,\\
\end{split}
\end{equation}
Then,  gathering these two  relations (\ref{e611082}) and (\ref{e611083}) we get
\begin{displaymath}
c^pB^{\mu\sigma}B^{-\mu^2}=B^{\sigma^2}B^{-\mu\sigma},
\end{displaymath}
thus
\begin{displaymath}
B^{(\sigma-\mu)^2}=c^p,\quad c\in K.
\end{displaymath}
In the other hand $B^{\sigma^{p-1}-1}=1$. But in the euclidean field ${\bf F}_p[X]$, we have $gcd((X^{p-1}-1), (X-\mu)^2))= X-\mu$.
Therefore $B^{\sigma-\mu}=\alpha_3^{ p},
\quad \alpha_3\in K$,
and so $\sigma(B)=B^{\mu}\times\alpha_3^{ p}$.
\end{enumerate}
The end of proof is similar to the
proof of previous lemma \ref{l203171} p. \pageref{l203171}.
\end{proof}
\end{thm}
%%% ====================================================================
%
%RRRRRRRR 260
%%% ====================================================================
\begin{thm}\label{t610241}
Let $B$ be the singular positive number defined in theorem \ref{l203272}.
If $m>\frac{p-1}{4}$ then  $B$ verifies the congruence $\modu p$:
\begin{equation}\label{e202211}
B^{p-1}\equiv 1-\delta\times
(\zeta+\mu^{-1}\zeta^u+\dots+\mu^{-(p-2)}\zeta^{u_{p-2}}) \modu \pi^{p-1},
\end{equation}
where $\delta\in \Z$ is coprime with $p$ when $B$ is non-primary
\begin{proof}
If $B$ is primary then it results of definition of primary numbers  . If $B$ is non-primary the proof is similar to theorem \ref{l202211} proof.
\end{proof}
\end{thm}
%%% ====================================================================
%
%RRRRRRRR 260
%%% ====================================================================
\section{On  $\pi$-adic expansion of singular units}
%%% ====================================================================
%
%RRRRRRRR 260
%%% ====================================================================
Let us fix $\eta$ for one of the  units $\eta_i$ of definition relation (\ref{e611085}).
The singular units $\eta$ verify  $\eta^{\sigma-\mu}\in O_{K^+}^{* p}$. Therefore
the results on  singular integers $B$ non-units  of section \ref{s108311} p. \pageref{s108311}
can be translated {\it mutatis mutandis} to get similar results  for  the  unit $p$-group
$U= O_{K^+}^*/O_{K^+}^{*p}$:
%%% ====================================================================
%
%RRRRRRRR 260
%%% ====================================================================
\begin{thm}\label{t610122}
Let $m$ be a natural number  $1\leq m\leq m$.
There exists singular units $\eta\in O_{K^+}^*$ verifying
$\sigma(\eta)=\eta^\mu\times \varepsilon^p$ with $\mu=u_{2m}$ and  $\varepsilon\in O_{K^+}^*$.
If $\eta$ is non-primary then $\pi^{2m} \ \|\ \eta^{p-1}-1$.
\end{thm}
%%% ====================================================================
%
%RRRRRRRR 260
%%% ====================================================================
\paragraph{The case $\mu=u_{2m}$ with  $ 2m>\frac{p-1}{2}$}\label{s201111}
\paragraph{ }
%%% ====================================================================
%
%RRRRRRRR 260
%%% ====================================================================
The next theorem for the $p$-unit group $U=O_{K^+}^*/O_{K^+}^{*p}$
is the translation of the similar theorem \ref{l202211} p. \pageref{l202211} for the singular negative integers.
\begin{thm}{ }\label{l203044}
Let $\mu=u_{2m},\quad p-3\geq 2m>\frac{p-1}{2}$.
The singular units  $\eta$  with $\eta^{\sigma-\mu}\in O_{K^+}^*$ verify the explicit congruence:
\begin{equation}\label{e203043}
\eta^{p-1}\equiv 1-\delta\times
(\zeta+\mu^{-1}\zeta^u+\dots+\mu^{-(p-2)}\zeta^{u_{p-2}}) \modu \pi^{p-1},
\end{equation}
where $\delta\in\Z$ is coprime with $p$ when $\eta$ is non-primary.
\end{thm}
%%% ====================================================================

%

%RRRRRRRR 260

%%% ====================================================================

Roland Qu\^eme

13 avenue du ch\^ateau d'eau

31490 Brax

France

mailto: roland.queme@wanadoo.fr

%$$$$$$$$$$$$$$$$$$$$$$$$$$$$
%$$$$$$$$$$$$$$$$$$$$$$$$$$$$
%$$$$$$$$$$$$$$$$$$$$$$$$$$$$

\end{document}